\documentclass[12pt, reqno]{amsart}
\UseRawInputEncoding
\usepackage{amsmath, amsthm, amscd, amsfonts, amssymb, graphicx, color}
\usepackage[bookmarksnumbered, colorlinks, plainpages]{hyperref}
\hypersetup{colorlinks=true,linkcolor=red, anchorcolor=green,
citecolor=cyan, urlcolor=red, filecolor=magenta, pdftoolbar=true}

\newcommand{\vf}{\varphi}
\newcommand{\A}{\mathfrak{A}}
\newcommand{\B}{\mathfrak{B}}
\newcommand{\M}{\mathfrak{M}}
\newcommand{\N}{\mathfrak{N}}
\newcommand{\h}{\mathcal{H}}
\newcommand{\cs}{\ensuremath{C^\ast}}
\newcommand{\lho}{\mathcal{B}_{A^{1/2}}(\h)}
\newcommand{\sa}{\mathcal{S}(\A)}
\newcommand{\sm}{\mathcal{S}(\M)}
\newcommand{\RR}{\mathcal{R}}
\newcommand{\dist}{{\rm dist}}
\newcommand{\NN}{\mathcal{N}}

\newcommand{\uno}{{\bf\rm \textbf{1}}}
\newcommand{\lh}{\mathcal{B}(\mathcal{H})}
\newtheorem{theorem}{Theorem}[section]
\newtheorem{lemma}[theorem]{Lemma}

\newtheorem{corollary}[theorem]{Corollary}
\theoremstyle{definition}
\newtheorem{definition}[theorem]{Definition}
\newtheorem{example}[theorem]{Example}

\theoremstyle{remark}
\newtheorem{remark}[theorem]{Remark}
\numberwithin{equation}{section}
\usepackage{color}
\newcommand{\cl}{{\rm cl}}
\newcommand{\gr}[1]{{\bf{#1}}}
\newcommand{\lhg}{\mathcal{B}(\gr R(A^{1/2}))}
\begin{document}

\setcounter{page}{1}
\title[$A$-spectral permanence property for \cs-algebras]
{$A$-spectral permanence property for \cs-algebras}

\author[M.~Mabrouk \MakeLowercase{and} A.~Zamani]
{Mohamed Mabrouk$^{1}$ \MakeLowercase{and} Ali Zamani$^{2}$}

\address{$^1$Faculty of Sciences of Sfax, Department of Mathematics University of Sfax Tunisia
\&
Department of Mathematics, Faculty of Applied Sciences, Umm Al-Qura University 21955 Makkah, Saudi Arabia
}
\email{msmabrouk@uqu.edu.sa}

\address{$^2$School of Mathematics and Computer Sciences, Damghan University, P.O.BOX 36715-364, Damghan, Iran}
\email{zamani.ali85@yahoo.com}

\begin{abstract}
For a positive element $A$ of a \cs-algebra $\A$, let ${\|X\|}_{A}$ denote the $A$-operator semi-norm of $X\in\A$.
In this paper, we aim to introduce and study the notion of $A$-spectrum for $X$ such that ${\|X\|}_{A}<\infty$.
In particular when $A$ is well-supported, we establish an $A$-spectral permanence property for \cs-algebras.
\end{abstract}
\maketitle
\section{Introduction and preliminaries}
Throughout this paper, let $\mathfrak{A}$ be a unital $C^*$-algebra with unit denoted by $\textbf{1}$
and we assume that $A\in\mathfrak{A}$ is a positive element.
Let $\mathfrak{A}^*$ and $\mathfrak{A}^+$ be the topological dual space of $\mathfrak{A}$
and the cone of positive elements in $\mathfrak{A}$, respectively.
We denote by $\lh$ the $C^*$-algebra of all bounded linear operators on a complex Hilbert
space $\big(\mathcal{H}, \langle \cdot, \cdot\rangle\big)$.
Let $\mathcal{R}(X)$ and $\mathcal{N}(X)$ stand for the range and the null space of $X\in\lh$, respectively.
By the Gelfand--Naimark theorem, every $C^*$-algebra $\mathfrak{A}$ can be considered as a norm
closed $*$-subalgebra of $\lh$ for some Hilbert space $\mathcal{H}$.
A $*$-subalgebra $\mathfrak{M}$ of $\lh$ is said to be a von Neumann algebra in $\mathcal{H}$ provided that
it is closed in the weak operator topology and contains the identity operator $\textbf{1}$.
For a fuller treatment of the theory of $C^*$-algebras and von Neumann algebras the reader is referred to \cite{Mu}.

A linear functional $f\in\mathfrak{A}^*$ is said to be positive, if $f(X)\geq 0$ for all $X\in\mathfrak{A}^+$.
Given a positive functional $f$ on $\mathfrak{A}$, we have the following useful version of the Cauchy--Schwarz inequality:
\begin{align*}
\big|f(X^*Y)\big|^2\leq f(X^*X)f(Y^*Y) \qquad (X, Y\in\mathfrak{A}).
\end{align*}
Let $\mathcal{S}(\mathfrak{A})$ denotes the set of states on $\mathfrak{A}$,
which is the set of all positive linear functionals $f$ on $\mathfrak{A}$ such that $\|f\|=1$.
Set $\mathcal{S}_{A}(\mathfrak{A})=\left\{\frac{f}{f(A)}: f\in\mathcal{S}(\mathfrak{A}),\, f(A)\neq0\right\}$
and for an element $X\in\mathfrak{A}$, let
${\|X\|}_{A}=\sup \left\{\sqrt{f(X^*AX)}: f\in\mathcal{S}_{A}(\mathfrak{A})\right\}$.
It is worth observing that ${\|\!\cdot\!\|}_{\textbf{1}}=\|\!\cdot\!\|$ and ${\|X\|}_{A}=0$ if and only if $AX=0$.
Further the set $\mathcal{S}_{A}(\mathfrak{A})$ is a non empty and convex subset of $\mathfrak{A}^*$,
but it is weak$^*$-compact if and only if $A$ is invertible in $\mathfrak{A}$;
see \cite[proposition~2.3]{B.M.Positivity}. In particular if $A$ is not invertible
and due to the lack of compactness of $\mathcal{S}_{A}(\mathfrak{A})$,
it may happen that ${\|X\|}_{A}=\infty$ for some $X\in\mathfrak{A}$; see \cite[Example~3.2]{B.M.Positivity}.
In the sequel  we will denote
$\mathfrak{A}^{A}=\left\{X\in\mathfrak{A}: \, {\|X\|}_{A} <\infty \right\}$.
The set $\mathfrak{A}^{A}$ is a subalgebra of $\mathfrak{A}$ not necessarily closed.
Also by \cite[Proposition~3.3]{B.M.Positivity} ${\|\!\cdot\!\|}_{A}$ is a semi-norm on $\mathfrak{A}^{A}$ and satisfies
${\|XY\|}_{A} \leq {\|X\|}_{A}{\|Y\|}_{A}$ for all $X, Y \in\mathfrak{A}^{A}$.
For an element $X\in\mathfrak{A}$, an element $Y\in\mathfrak{A}$ is said to be an $A$-adjoint of $X$ if $AX=Y^*A$.
An element $X\in\mathfrak{A}$ is said to be $A$-self-adjoint if $AX$ is self-adjoint, i.e., $AX = X^*A$.
Denote by $\mathfrak{A}_{A}$ the set of all elements in $\mathfrak{A}$ that admit $A$-adjoints.
It is well known that $\mathfrak{A}^A$ is a subalgebra of $\A$ which is neither closed nor dense in $\A$ even when $\A=\lh$.
Also, $\A_A=\A$ if $A$ is in the center of $\A$ but in general $\A_A\neq \A$
and that $\A_A\subseteq \A^A$ and $\A^A=\lho$ when $\A=\lh$.
Here $\lho$ is the set of all bounded linear operators on $\mathcal{H}$ having an $A^{1/2}$-adjoint.
More properties of $\mathfrak{A}_{A}$ were investigated in \cite{amz, B.M.Positivity}.

Recently in \cite{Baklouti.Namouri} and for $X\in\lho$, Baklouti and Namouri defined
the notion of $A$-spectrum $\sigma_A(X, \lh)$ (see also \cite{Feki1}).
They established some permanence properties of this concept.
The current paper is motivated by \cite{Baklouti.Namouri}.
We shall study the notion of $A$-spectrum $\sigma_A(X,\A)$ for an element $X\in\A^A$.

The contents of the paper are the following. Section \ref{sec1} contains basic results on
$A$-spectrum in \cs-algebras. Several properties are obtained. In Section
\ref{sec2} and if  $A$ is well-supported we prove an $A$-spectral permanence theorem for \cs-algebras.
Namely, we show that if $\A$ and $\B$ are two von Neumann algebras such that $\B\subseteq \A$ and $A\in\B$
then $\sigma_A(X,\A)=\sigma_A(X,\B)$ for any $X\in \B^A$. Finally, some corollaries are presented.
\section{$A$-invertibility on $\A^A$}\label{sec1}
We shall denote \cs-algebras by $\A$ and $\B$ and von Neumann algebras by $\M, \N\cdots$.
Recall that, given an element $X$ in $\A$ we say that $Y\in\A$ is an $A$-adjoint of $X$ if $Y$
satisfies the equation $AX = Z^*A$.
In the frame of von Neumann algebras, this kind of equations can
be studied by applying the next theorem (called Douglas theorem for von Neumann algebras).
\begin{lemma}\cite{Nayak.2}\label{lemDouglas}
Let $X, Y \in \M$.
Then the following conditions are equivalent:
\begin{itemize}
\item[(i)] $X^*X\leq \alpha Y^*Y$ for some $\alpha\geq0$;
\item[(ii)] there exists $Z\in \M$ such that $X=Z^*Y$.
\end{itemize}
Moreover, if $X^*X=Y^*Y$, then $Z$ can be chosen to be a partial isometry with initial
projection the range projection of $X$, and final projection as the range projection of $Y$.
\end{lemma}
Note that Lemma \ref{lemDouglas} fails to be true in the general setting of \cs-algebras.
In fact, for a locally compact Hausdorff space $\Omega$,
the commutative \cs-algebra $C_0(\Omega)$ has the Douglas property if and only if $\Omega$
is sub-Stonean. For more details the reader may consult \cite{Nayak.1} and \cite{Nayak.2}.
Nevertheless, we have the following.
\begin{lemma}\label{vncs}
Let $X\in\A$ and $B\in\mathfrak{A}^+$. If $X^*X\leq B$ and $ 0 < \alpha < 1/2$,
then there is a $V\in \cl(\A A)\subseteq\A$ with $V^*V\leq B^{1-2\alpha}, VV^*\leq (XX^*)^{1-2\alpha}$ and $X = VB^\alpha$.
\end{lemma}
\begin{proof}
See \cite[Proposition II.3.2.1]{blackadar2006operator}.
\end{proof}
Recall that an element $X\in\A$ is said to be well-supported if $\sigma(X^*X)\backslash\{0\}$ is closed.
See \cite[Page~69]{blackadar2006operator} for more information.
As mentioned in \cite{Ar.Co.Go}, by Douglas theorem and when $\A=\lh$, we have $\A^A=\lho$.
In particular if $A$ is well-supported we see that
$\A_A=\A^A$.
But in general $\A_A\subseteq \A^A\subseteq \lho$.
The next theorem extends \cite[Proposition~1.1]{Ar.Co.Go}
and gives another characterization of the set $\A^A$.
\begin{theorem}\label{adjoint}
The following statements hold.
\begin{itemize}
\item[(i)] $\A^A=\Big\{X\in \A:  \ \exists\ c>0 , \ U\in\A\ \text{such that} \  A^{1/2}X=UA^{1/4}\ \text{and}\  U^*U\leq cA^{1/2}\Big\}$.
\item[(ii)] If $A$ is well-supported then
\begin{align*}
\A^A=\lho\cap\A =\Big\{X\in \A: X \ \text{has an} \ A^{1/2}-\text{adjoint}\Big\}.
\end{align*}
\item[(iii)] $\M^A=\lho\cap\M =\Big\{X\in \M: \ X \ \text{has an} \ A^{1/2}-\text{adjoint}\Big\}$.
\end{itemize}
\end{theorem}
\begin{proof}
(i) If $X\in \A$ is such that $ A^{1/2}X=UA^{1/4}$ for some $U\in\A$ and $c>0$ with $U^*U\leq cA^{1/2}$.
Then for any $f\in\sa$, since $A^{1/4}U^*UA^{1/4}\leq cA$, we have
\begin{align*}
f(X^*AX)=f(A^{1/4}U^*UA^{1/4})\leq c f(A),
\end{align*}
and so ${\|X\|}_{A}<\infty$. Hence $X\in\A^A$.
Now, let $X\in\A^A$ and set $c={\|X\|}^2_{A}$. Then by the definition of ${\|X\|}_{A}$ we have
$f(X^*AX)\leq c f(A)$ for any state $f\in\sa$ such that $f(A)\neq 0$.
Now, let $f\in\sa$ such that $f(A)=0$. As $A\neq 0$, then there exists $g\in\sa$ so that $g(A)\neq 0$.
Set $f_{\varepsilon}=\varepsilon f+(1-\varepsilon)g$, then $f_{\varepsilon}$ is a state on $\A$ and
$f_{\varepsilon}(A)\neq 0$ for any $0<\varepsilon<1$. In particular  $f_{\varepsilon}(X^*AX)\leq c f_{\varepsilon}(A)$
for any $0<\varepsilon<1$. Tending $\epsilon$ to $1$ we get $f(X^*AX)\leq c f(A)$ for any $f\in\sa$.
Hence $X^*AX\leq cA$. Applying Lemma \ref{vncs} with $\alpha=\frac{1}{4}$ and $B=cA$
we get $A^{1/2}X=c^{1/4}VA^{1/4}$ for some $V\in \cl(\A A)\subseteq\A$ such that $V^*V\leq c^{1/2} A^{1/2}$.
Put $U=c^{1/4}V$, we see that $A^{1/2}X=UA^{1/4}$ and $U^*U\leq cA^{1/2}$.

(ii) Assume that $A$ is well-supported. Then $\cl(\A A)=\A A$, by \cite[Proposition. II.3.2.11]{blackadar2006operator}.
Hence $U= BA$ for some $B\in\A$, by (i). Therefore
\begin{align*}
A^{1/2} X=BA^{5/4}=\left(A^{3/4}B^*\right)^*A^{1/2}.
\end{align*}
The proof is thus complete since $A^{3/4}B^*$ is an $A^{1/2}$-adjoint of $X$.

(iii) The proof is so similar to that of (ii) that we omit it.
\end{proof}
Let $A$ be well-supported. Recall that the Moore--Penrose inverse of $A$, denoted by $A^\dag$, is the unique element $Z\in\A$ that satisfies
\begin{align*}
A = AZA, \, ZAZ = Z, \,(AZ)^* = AZ, \, (ZA)^*= ZA.
\end{align*}
Note that if we denote by $P$ the orthogonal projection onto the range of $A$, then $AA^\dag=A^\dag A=P$.
See \cite{Roch.Silbermann} for more information.
\begin{corollary}\label{corn}
If $A$ is well-supported then $\A^A=\A_A$.
\end{corollary}
\begin{proof}
It suffices to show that $\A^A\subseteq\A_A$. To do so let $X\in\A^A$.
By Theorem \ref{adjoint}(ii) $A^{1/2}X=Y^*A^{1/2}$ for some $Y\in \A^A$.
Keeping in mind that $A^\dag\in\A$ and $AA^\dag=A^\dag A=P$, we get
\begin{align*}
AX=A^{1/2}Y^*A^{1/2}=\left(A^{1/2}Y^*A^{1/2}A^\dag\right)A = \left((A^\dag)^*A^{1/2}YA^{1/2}\right)^*A.
\end{align*}
Hence $(A^\dag)^*A^{1/2}YA^{1/2}$ is an $A$-adjoint of $X$ and so $X\in\A_A$ as desired.
\end{proof}
\begin{remark}\label{rem1}
Observe that by the proof of Theorem \ref{adjoint} we have $f(X^*AX)=0$ for any $X\in\A^A$ and $f\in\mathcal{S}(\A)$ whenever $f(A)=0$.
\end{remark}
Recently, the authors in \cite{Baklouti.Namouri} introduced the concept $A$-invertibility of elements in $\lho$.
More precisely, an operator $X\in\lho$ is said to be
$A$-invertible in $\lho$ if there exists $Y\in\lho$ such that $AXY=AYX=A$.
A necessary and sufficient condition for the existence of $A$-inverses in $\lho$ is stated in \cite[Theorem~4.2]{Baklouti.Namouri} as follows.
\begin{lemma}\cite{Baklouti.Namouri}\label{hs}
An operator $X\in\lho$ is invertible in $\lho$ if and only if the following two conditions are satisfied:
\begin{itemize}
\item[(i)] there exists $c > 0$ such that $\frac{1}{c}{\|h\|}_{A}\leq {\|Xh\|}_{A} \leq {\|h\|}_{A}$ for all $h\in\h$;
\item[(ii)] $\RR(AX) = \RR(A)$.
\end{itemize}
\end{lemma}
Now, we introduce the following definition.
\begin{definition}
We say that an element $X\in\A^A$ is $A$-left (resp. $A$-right) invertible in $\A^A$
if there exists $Y\in\A^A$ such that $AYX=A$ (resp. $AXY=A$).
An element $X\in\A^A$ is said to be $A$-invertible in $\A^A$ if it has an $A$-inverse in $\A^A$, that
is, if there exists an element $Y\in\A^A$ such that $AXY=AYX=A$.
\end{definition}
\begin{remark}
Observe that a given $X\in\A^A$ may admit none, one or many $A$-inverses.
In fact, if $Y$ is an $A$-inverse of $X$ and $AXZ =AZX= 0$ for
some $Z \in\A$ then $Y + Z$ is also an $A$-inverse of $X$.
\end{remark}
It is clear that any $A$-invertible element in $\A^A$ is $A$-left and $A$-right invertible in $\A^A$.
The converse is also true as the following result shows.
\begin{theorem}\label{LRinverses}
Let $X\in\A^A$. The following statement are equivalent.
\begin{enumerate}
\item[(i)] $X$ is $A$-invertible in $\A^A$;
\item[(ii)] $X$ is $A$-left and $A$-right invertible in $\A^A$;
\item[(iii)] there exists $Y\in\A^A$ such that $AXY=A$ and $\NN(AX)\subseteq\NN(A)$.
\end{enumerate}
\end{theorem}
\begin{proof}
(i)$\Longrightarrow$(ii) The implication is trivial.

(ii)$\Longrightarrow$(iii) Assume that $AXY=AZX=A$ for some $Y$ and $Z$ in $\A^A$.
Note that
\begin{align}\label{Z.I}
\NN(A^t)=\NN(A),
\end{align}
for any $t>0$ (see, for example, \cite[Lemma~2.9]{Z.LAMA}).
By Theorem \ref{adjoint}(i), there exists $V\in\A$ such that $A^{1/2}Z = VA^{1/4}$.
We have
\begin{align}\label{Z.Ii}
A=AZX=A^{1/2}\big(A^{1/2}Z\big)X=A^{1/2}VA^{1/4}X.
\end{align}
Now, let $h\in\NN(AX)$. Then $Xh\in\NN(A)$ and by \eqref{Z.I} it follows that $Xh\in\NN(A^{1/4})$.
Thus $h\in\NN(A^{1/4}X)$ and
therefore by \eqref{Z.Ii}, we get $h\in\NN(A)$. Hence $\NN(AX)\subseteq\NN(A)$.

(iii)$\Longrightarrow$(i) Suppose (iii) holds. Hence $AXY=A$ for some $Y\in\A^A$ and $\NN(AX)\subseteq\NN(A)$.
Let $P$ be the orthogonal projection onto the range of $A$. By Theorem \ref{adjoint}(i), there exists $U\in\A$ such that $A^{1/2}X = UA^{1/4}$.
Let $h\in\NN(A)$. So, by \eqref{Z.I}, $h\in\NN(A^{1/4})$. Hence $h\in\NN(UA^{1/4})=\NN(A^{1/2}X)$ and then $h\in\NN(AX)$.
Therefore $\NN(AX)=\NN(A)$. Similarly we can see that $\NN(AY)=\NN(A)$.
We have $h\in\NN(AYX)$ if and only if $Xh\in\NN(AY)$, or equivalently, if and only if $Xh\in\NN(A)$.
This occurs if and only if $h\in\NN(AX)$, or equivalently, $h\in\NN(A)$.
Whence $\NN(AYX)=\NN(A)$. Keeping in mind that $X\big(\NN(A)\big)\subseteq \NN(A)$ and $Y\big(\NN(A)\big)\subseteq \NN(A)$,
we get $(PYX)^2 = PYX$. That is $PYX$ is a projection on $\h$.
Finally, since $\RR(PYX)\subseteq\RR(P)=\cl(\RR(A))$ and $\NN(AYX)=\NN(A)$, it follows that $PYX=P$. We have
\begin{align*}
AYX=APYX=AP=A
\end{align*}
and so $AYX=A$. Therefore, $X$ is $A$-invertible.
\end{proof}
Since a von Neumann algebra contains the range projections of their elements, we also get the following result.
\begin{corollary}
Let $X\in\M^A$ and $P$ be the orthogonal projection onto the range of $A$.
The following conditions are equivalent:
\begin{itemize}
\item[(i)] $XP$ is $A$-invertible;
\item[(ii)] $PX$ is $A$-invertible;
\item[(iii)] $X$ is $A$-invertible.
\end{itemize}
\end{corollary}
\begin{proof}
(i)$\Longrightarrow$(ii) Suppose (i) holds. Then $AXPY=AYXP=A$ for some $Y\in\M^A$.
Since $A=AP$, $PXP=PX$ and $PYP=PY$ then, $A=A(PX)Y$ and $A(PY)(PX)=A$. Whence $PX$ is $A$-invertible, by Theorem \ref{LRinverses}.

(ii)$\Longrightarrow$(iii) Let $PX$ be $A$-invertible. Then there exists $Y\in\A^A$ such that $APXY=AYXP=A$.
Keeping in mind that $A=AP$ and $PYXP=PYX$ it yields that $AXY=A$ and $AYX = APYX = A(PYXP) = AYXP=A$.
Whence $X$ is $A$-invertible.

The implication (iii)$\Longrightarrow$(i) can be handled in a similar way.
\end{proof}

Let us quote a result from \cite{FILLMORE1971254}.
\begin{lemma}\cite[Corollary 1]{FILLMORE1971254}\label{FILLMORE}
Let $B$ and $C$ be bounded operators on $\h$. There
exists an invertible operator $D$ on $\h$ such that $B = CD$ if and only if $B$ and $C$ have the same range and nullity.
\end{lemma}
The next result shows that the $A$-inverse in $\lho$ may be chosen to be invertible in $\lh$.
\begin{theorem}\label{inv6}
Let $X\in\lho$. The following statement are equivalent.
\begin{enumerate}
\item[(i)] $X$ is invertible in $\lho$;
\item[(ii)] there exists an operator $Y\in\lho$ which is invertible in $\lh$ such that $AXY =AYX =A$.
\end{enumerate}
\end{theorem}
\begin{proof}
The implication (ii)$\Longrightarrow$(i) is trivial. Let us show the implication (i)$\Longrightarrow$(ii).
Assume that $X$ is invertible in $\lho$. So, by Lemma \ref{hs},
$\RR(AX) = \RR(A)$ and there exists $c > 0$ such that
\begin{align}\label{inv6.1}
\frac{1}{c}{\|h\|}_{A}\leq {\|Xh\|}_{A} \leq {\|h\|}_{A} \qquad (h\in\h).
\end{align}
From this it follows that $\NN(AX)=\NN(A)$. Hence, by Lemma \ref{FILLMORE}
there exists an invertible operator $Y$ on $\h$ such that $AXY=A$.
Now, let $k\in\h$. By \eqref{inv6.1} we have
\begin{align*}
{\|Yk\|}^2_{A}&\leq c{\|XYk\|}^2_{A} = c\langle AXYk, XYk\rangle
\\&= c\langle Ak, XYk\rangle = c\langle k, AXYk\rangle = c\langle k, Ak\rangle = c{\|k\|}^2_{A},
\end{align*}
and so ${\|Yk\|}_{A}\leq \sqrt{c}{\|k\|}_{A}$. From this it follows that $Y\in\lho$.
Hence we have $AXY=A$ with $Y\in\lho$ and $\NN(AX)\subseteq\NN(A)$. Now, Theorem \ref{LRinverses} leads up to (ii).
\end{proof}
Now, by Theorem \ref{LRinverses}, the definition of $A$-spectrum for an element $X\in\A^A$ is naturally.
More precisely set:
\begin{align*}
\sigma_A^l(X, \A)=\Big\{\lambda \in\mathbb{C}: \lambda \uno-X\ \text{is not}\ A-\text{left invertible}\Big\},
\end{align*}
\begin{align*}
\sigma_A^r(X, \A)=\Big\{\lambda \in\mathbb{C}: \lambda \uno-X\ \text{is not}\ A-\text{right invertible}\Big\},
\end{align*}
and
\begin{align*}
\sigma_A(X, \A)=\sigma_A^l(X, \A)\cup\sigma_A^r(X, \A).
\end{align*}
We call $\sigma_A^l(X, \A)$ and $\sigma_A^r(X, \A)$ and $\sigma_A(X, \A)$
the left $A$-spectrum and the right $A$-spectrum and the $A$-spectrum of $X\in\A^A$, respectively.
Also, if no ambiguity arises we shall simply denote these subsets simply by $\sigma_A^l(X)$ and $\sigma_A^r(X)$ and $\sigma_A(X)$.
\begin{remark}
It is worth noting that when $A$ is invertible in $\A$, the above sets coincide with the left spectrum $\sigma^l(X)$,
the right spectrum $\sigma^r(X)$ and the spectrum $\sigma(X)$ of $X$, respectively.
\end{remark}
The next theorem is a variant of \cite[Theorem~4.2]{Baklouti.Namouri}.
It characterizes $A$-invertibility on $\A^A$.
\begin{theorem}\label{thvn}
Let $\M$ be a von Neumann algebra on $\h$ (resp. a \cs-algebra with $A$ well-supported) and $X\in\M^A$ (resp. $X\in\A^A$).
Then $X$ is $A$-invertible in $\M^A$ (resp. in $\A^A$) if and only if the following two conditions hold:
\begin{itemize}
\item[(i)] there exists $c>0$ such that $\frac{1}{c}f(A)\leq f(X^*AX)\leq cf(A)$ for any $f\in\sm$ (resp. $f\in\sa$);
\item[(ii)] there exists $\alpha>0$ such that $A^2\leq \alpha AXX^*A$.
\end{itemize}
\end{theorem}
\begin{proof}
We shall give the proof only for von Neumann algebra.
The case of a \cs-algebra $\A$ with $A$ well-supported uses the same scheme together with Theorem \ref{adjoint}(ii).

Assume that $X$ is $A$-invertible in $\M^A$. Hence, $AXY=AYX=A$ for some $Y\in\M^A$.
Put $c:= \max\{{\|X\|}_{A}, {\|Y\|}_{A}\}$ and $\alpha:=\|Y\|^2$.
Since $X\in\M^A$, then $f(X^*AX)\leq {\|X\|}_{A}f(A)$ for any $f\in\sm$ such that $f(A)\neq 0$.
Recall that by Remark \ref{rem1}, we have $f(X^*AX)=0$ when $f(A)=0$,
and hence $f(X^*AX)\leq {\|X\|}_{A}f(A)$. Thus $f(X^*AX) \leq cf(A)$ for any $f\in\sm$.

Now, if $f(X^*AX)=0$ again Remark \ref{rem1} entails that $f(A)=0$ and so $\frac{1}{c}f(A)=f(X^*AX)$.
If $f(X^*AX)\neq0$, then by the Cauchy--Schwarz inequality we have
\begin{align*}
f^2(A)=f^2(AYX)&\leq f(A)f(X^*Y^*AYX)
\\&= f(A)f(X^*AX)\frac{f(X^*Y^*AYX)}{f(X^*AX)}
\\& \leq f(A)f(X^*AX){\|Y\|}_{A} \leq cf(A)f(X^*AX),
\end{align*}
whence $f^2(A)\leq cf(A)f(X^*AX)$. Therefore, $\frac{1}{c}f(A)\leq f(X^*AX)$ for any $f\in\ \sm $.

Also, since $YY^*\leq {\|Y\|}^2\uno = \alpha \uno$, by \cite[Theorem~2.2.5]{Mu}, we have
\begin{align*}
AX\big(\alpha \uno-YY^*\big)X^*A\geq 0.
\end{align*}
Thus,
\begin{align*}
A^2 &= AXYY^*X^*A
\\& = \alpha AXX^*A - AX\big(\alpha \uno-YY^*\big)X^*A \leq \alpha AXX^*A,
\end{align*}
and so $A^2\leq \alpha AXX^*A$.

For the converse, note that condition (i) implies that
\begin{equation}\label{wh}
\frac{1}{c}\langle Ah, h\rangle \leq \langle X^*AX h, h\rangle \leq c \langle Ah, h\rangle \quad (h\in\h).
\end{equation}
Since $A^2\leq \alpha AXX^*A$ for some $\alpha > 0$, then by Lemma \ref{lemDouglas} there exists $Y\in\M$ such
that $A=Y^*X^*A$ or, equivalently $AXY=A$. Using the condition \eqref{wh}, it is clear that $\NN(AX)\subseteq\NN(A)$.
Also, for any $h\in\h$ by \eqref{wh}, we have
\begin{align*}
\langle Y^*AYh, h\rangle&=\langle AYh, Yh\rangle
\\& \leq c\langle X^*AXYh, Yh\rangle
\\&=c\langle Y^*X^*AXY h, h\rangle = c\langle AXYh, h\rangle
=c\langle Ah, h\rangle.
\end{align*}
Hence ${\|Y\|}_A<\infty$ and then $Y\in\M^A$.
This together with Theorem \ref{LRinverses} allows us to conclude.
\end{proof}
The next result investigates the $A$-invertibility of the product of two $A$-invertible elements in a von Neumann algebras.
\begin{theorem}\label{invprof}
Let $X$ and $Y$ in $\M^A$ be $A$-invertible with $A$-inverses $Z\in\M^A$
and $W\in\M^A$, respectively. Then $XY$ is $A$-invertible with an $A$-inverse $WZ$.
\end{theorem}
\begin{proof}
Assume that $AXZ=AZX=A$ and $AYW=AWY=A$. By Theorem \ref{LRinverses} we have $\NN(AX)\subseteq\NN(A)$ and $\NN(AY)\subseteq\NN(A)$.
Let $h\in\NN(AXY)$.
Hence $Yh\in\NN(AX)$ and then $Yh\in\NN(A)$. Thus $h\in\NN(AY)$ and so $h\in\NN(A)$. Therefore, $\NN(AXY)\subseteq\NN(A)$.
Further, by \eqref{Z.I}, we have $A^{1/2}XZ=A^{1/2}ZX=A^{1/2}$ and $A^{1/2}YW=A^{1/2}WY=A^{1/2}$.
Since, by Theorem \ref{adjoint}(iii), we have $A^{1/2}X=L^*A^{1/2}$ for some $L\in\M^A$, then
\begin{align*}
AXYWZ=A^{1/2}L^*A^{1/2}YWZ=A^{1/2}L^*A^{1/2}Z=AXZ=A.
\end{align*}
This completes the proof, since $WZ\in\M^A$.
\end{proof}
We close this section by the following result.
\begin{theorem}
If $X, Y, R$ and $W$ are in $\A^A$ such that $A^{1/2}X=W^*A^{1/2}$ and $A^{1/2}Y=R^*A^{1/2}$,
then $X$ is $A$-invertible with $A$-inverse $Y$ if and only if $W$ is $A$-invertible with $A$-inverse $R$.
\end{theorem}
\begin{proof}
Assume that $AXY=AYX=A$. Hence, by \eqref{Z.I}, we have
$A^{1/2}XY=A^{1/2}YX=A^{1/2}$ and then $W^*R^*A^{1/2}=R^*W^*A^{1/2}=A^{1/2}$.
Therefore, $A^{1/2}RW=A^{1/2}WR=A^{1/2}$ and so $ARW=AWR=A$. Whence $W$ is $A$-invertible with $A$-inverse $R$.
The converse can be handled similarly.
\end{proof}
\section{$A$-spectral permanence}\label{sec2}
One of the most important results in the theory of \cs-algebra is the so called spectral permanence theorem (see \cite[Theorem. 2.1.11]{Mu}),
which states that if $\B$ is a \cs-subalgebra of a unital \cs-algebra $\A$ containing the unit of $\A$,
then $\sigma(X, \B) = \sigma(X, \A) \ (X\in\B)$.
This means that if $X\in\B$ is invertible in $\A$, then $X$ is invertible in $\B$. Evidently, if $A$ is invertible,
then $\sigma_A(X)=\sigma(X)$ for any $X\in\A$ and the spectral permanence property holds.
But in general, this property fails to be true for the $A$-spectrum as the following example shows.
\begin{example}\label{e.009}
Let $\Omega := \big\{\frac{1}{n}: n\ge 1\big\}\cup\{0\}$. Note that $\Omega$ is a compact
Hausdorff space in the topology inherited from $\mathbb{R}$ with the usual topology.
Let $A$ and $X$ be the two functions in the
commutative \cs-algebra $C(\Omega)$ defined as follows:
\begin{align}\label{nws}
A(t)=\begin{cases}
t & \text{if}\ t=\frac{1}{2n}, n\ge 1,\\
0 & \text{otherwise},
\end{cases}
\end{align}
and $X(t)=t$ for any $t\in \Omega$. If $Y$ is a complex-valued function such that $A=AXY$,
then we must have $Y(\frac{1}{2n})=2n$ for any $n\ge 1$.
Hence $Y$ cannot be continuous.
On the other hand the function
\begin{align*}
Y(t)=\begin{cases}
\frac{1}{t} & \text{if}\ t=\frac{1}{2n}, n\ge 1,\\
0 & \text{otherwise }
\end{cases}
\end{align*}
is an $A$-inverse of $X$ in $L^\infty(\Omega, \mu)$, where $\mu$ is the Lebesgue measure.
Hence $0\notin\sigma_A\left(X, C(\Omega)\right)$ but $0\in\sigma_A\left(X, L^\infty(\Omega, \mu)\right)$.
\end{example}
Observe that the positive operator given by \eqref{nws} is not well-supported since $0$ is an accumulation point of $A(\Omega)$.
In what follows our purpose is to prove an $A$-spectral permanence theorem for \cs-algebras when $A$ is well-supported.
Whence, in the remainder of the paper, we suppose that $A$ is well-supported.
En route to prove our promised result, we begin with some preparatory results needed for our investigation.
The first lemma is given in \cite[Theorem ~2.13]{amz}.
\begin{lemma}\label{fromAMZ}
Let $g: A\A^A\longrightarrow \mathbb{C}$ such that $g(A)=1$ and $|g(AY)|\leq {\|Y\|}_{A}$ for all $Y\in\A^A$.
Then there exists $f\in\mathcal{S}_{A}(\A)$ such that $f(AY)=g(AY)$ for all $Y\in\A^A$.
\end{lemma}
The second lemma reads as follows.
\begin{lemma}\label{norme1}
Let $\A\subseteq\lh$ and $X\in\A^A$.
If ${\|X\|}_{A}<1$, then $\uno-X$ is $A$-invertible in $\A^A$.
\end{lemma}
\begin{proof}
Note that $A^{1/2}$ is well-supported.
Set $W=(A^{1/2})^\dag X^* A^{1/2}$. Since $X(\NN(A))\subset \NN(A)$, we have $A^{1/2}X = W^*A^{1/2}$.
Thus $W$ is an $A^{1/2}$-adjoint of $X$.
Further, by \cite[Theorem 3.5]{B.M.Positivity} and \cite[Proposition 2.2]{arias2008metric},
we have ${\|X\|}_{A}=\|A^{1/2}X(A^{1/2})^\dag\|=\|W\|$.
Applying \cite[Theorem 1-2.2]{Mu}, we see $ \uno-W$ is invertible in $\A$ with
\begin{equation*}
\left(\uno-W\right)^{-1}=\sum_{n=0}^{\infty}W^n=\sum_{n=0}^{\infty}(A^{1/2})^\dag(X^*)^n A^{1/2}.
\end{equation*}
It is clear that $(\uno-W)^{-1}\in\A^A$ and whence $\uno-W$ is $A$-invertible.
Finally, by using the fact that $A^{1/2}(\uno-W)=(\uno-X)^*A^{1/2}$ one
can see easily that $\uno-X$ is $A$-invertible in $\A^A$.
\end{proof}
\begin{corollary}
Let $(X_n)_{n\geq 1}$ be a sequence of $A$-invertible elements of $\A^A$
and $X\in\A^A$ be a non $A$-invertible element such that $\displaystyle\lim_{n\to\infty}{\|X-X_n\|}_{A}=0$.
If $Y_n\in\A^A$ is an $A$-inverse of $X_n$ for each $n\geq 1$, then $\displaystyle\lim_{n\to\infty}{\|Y_n\|}_{A}=+\infty$.
\end{corollary}
\begin{proof}
Suppose to the contrary that $\displaystyle\lim_{n\to\infty}{\|Y_n\|}_{A}<\infty$.
Then there exist $c> 0$ and a subsequence $(Y_{\theta(n)})_{n\geq 1}$ of $(Y_n)_{n\geq 1}$ such that
${\|Y_{\theta(n)}\|}_{A} \leq c$ for any $n\geq 1$.
We have
\begin{align*}
AX=A\left(X-X_{\theta(n)}\right)+AX_{\theta(n)}= AX_{\theta(n)}\Big(Y_{\theta(n)}\left(X-X_{\theta(n)}\right)+\uno\Big).
\end{align*}
Since ${\left\|Y_{\theta(n)}\left(X-X_{\theta(n)}\right)\right\|}_{A}\leq c{\left\|X-X_{\theta(n)}\right\|}_{A}$ and
$\displaystyle\lim_{n\to\infty}{\left\|X-X_{\theta(n)}\right\|}_{A}=0$, then, by Lemma \ref{norme1},
$Y_{\theta(n)}\left(X-X_{\theta(n)}\right)$ is $A$-invertible for $n$ large enough.
Hence, by Theorem \ref{invprof} $X$ is $A$-invertible and we have a contradiction.
\end{proof}
The third lemma reads as follows.
\begin{lemma}\label{lem-1}
Let $X\in\A^A$. The following statements hold.
\begin{itemize}
\item[(i)] If $AX\A^A\varsubsetneq A\A^A$, then there exists $f\in\mathcal{S}_{A}(\A)$ such that $f(AXY)=0$ for all $Y\in\A^A$.
\item[(ii)] If $A\A^AX\varsubsetneq A\A^A$, then there exists $f\in\mathcal{S}_{A}(\A)$ such that $f(AYX)=0$ for all $Y\in\A^A$.
\end{itemize}
\end{lemma}
\begin{proof}
(i) Put $N(AZ):={\|Z\|}_{A}$ for any $Z\in\A^A$ and note that $N(\cdot)$ defines a norm on $A\A^A$. We have
\begin{align*}
\dist(A, AX\A^A)=\displaystyle{\inf_{Z\in\A^A}N(AXY-A)}=\displaystyle{\inf_{Z\in\A^A}{\|XZ-\uno\|}_{A}}.
\end{align*}
It is clear that $\dist(A, AX\A^A)\leq {\|\uno\|}_{A}=1$. If $\dist(A, AX\A^A)<1$, then ${\|XZ-\uno\|}_{A}<1$ for some $Z\in\A^A$.
By Lemma \ref{norme1}, $XZ$ is $A$-invertible and then $AX\A^A=A\A^A$ which contradicts our assumption.
Therefore $\dist(A, AX\A^A)=1$. Hence there exists a bounded linear functional $g$ of norm $1$
on $A\A^A$ such that $g_{|AX\A^A}=0$ and $g(A) = 1$.
Since $|g(AZ)|\leq N(AZ)={\|Z\|}_A$ for any $Z\in\A^A$,
then by Lemma \ref{fromAMZ} there exists $f\in\mathcal{S}_{A}(\A)$ such that $f(AZ) = g(AZ)$ for all $Z\in\A^A$.
Thus $f(AXY)= g(AXY)=0$ for all $Y\in\A^A$.

(ii) The proof of (ii) is similar (i) and so we omit it.
\end{proof}
Now, we are in a position to state one of the main results of this section.
\begin{theorem}\label{th2-2}
Let $X\in\A^A$. Then
\begin{align}\label{sup1}
\sigma_A^r(X)=\Big\{f(AX): f\in\mathcal{S}_{A}(\A)\ \text{and}\  f(AXY)=f(AX)f(AY), \forall Y\in\A^A\Big\}
\end{align}
and
\begin{align}\label{sup2}
\sigma_A^l(X)=\Big\{f(AX): f\in\mathcal{S}_{A}(\A)\ \text{and}\  f(AYX)=f(AX)f(AY), \forall Y\in\A^A\Big\}.
\end{align}
\end{theorem}
\begin{proof}
First, let $\lambda\in\sigma_A^r(X)$. Then $A\left(\lambda\uno-X\right)\A^A\subsetneq A\A^A$.
By Lemma \ref{lem-1}(i), there exists $f\in\mathcal{S}_{A}(\A)$ such that $f(A(\lambda\uno-X)Y)=0$ for all $Y\in\A^A$.
Accordingly, $f(AXY)=\lambda f(AY)$ for all $Y\in\A^A$. Since $f(A)=1$ we get by taking $Y=\uno\in\A^A$,
$\lambda=f(AX)$ and then $f(AXY)=f(AX)f(AY)$ for every $Y\in\A^A$.
Conversely, if $\lambda=f(AX)$ for some $f\in\mathcal{S}_{A}(\A)$ such that $f(AXY)=f(AX)f(AY)$ for all $Y\in\A^A$.
Then $X-\lambda\uno$ is not $A$-invertible since if $A(X-\lambda\uno)Y=A$ for some $Y\in\A^A$,
we get $f(A)=f\left(A(X-\lambda\uno)Y\right)=f(AX)f(AY)-\lambda f(AY)=0$ which is impossible since $f(A)=1$.
In a similar way one can prove the equality for $\sigma_A^l(X)$.
\end{proof}
\begin{remark}
Observe that if $A$ is not well-supported, then
\begin{align*}
\Big\{f(AX): f\in\mathcal{S}_{A}(\A)\ \text{and}\  f(AXY)=f(AX)f(AY), \forall Y\in\A^A\Big\}\varsubsetneq \sigma_A^r(X)
\end{align*}
or
\begin{align*}
\Big\{f(AX): f\in\mathcal{S}_{A}(\A)\ \text{and}\  f(AYX)=f(AX)f(AY), \forall Y\in\A^A\Big\}\varsubsetneq \sigma_A^l(X).
\end{align*}
In fact the equalities hold if and only if $A$ is well-supported.
Indeed: if $A$ is well-supported then the result follows from Theorem \ref{th2-2}.
For the converse, assume that equalities \eqref{sup1} and \eqref{sup2} are true for any $X$.
In particular $0\notin\sigma_A(A)$. That is $A$ is $A$-invertible. Whence there exists $B\in\A^A$ so that $A=A^2B=ABA$.
So $A\A$ is closed and then by \cite[Theorem~2.1]{Kulkarni2000ACO} $A$ is well-supported.
\end{remark}
When $A=\uno$ we get the following result.
\begin{corollary}\label{cor-2}
Let $X\in\A$. Then
\begingroup\makeatletter\def\f@size{9}\check@mathfonts
\begin{equation*}
\sigma(X)=\Big\{f(X): f\in\mathcal{S}(\A)\ \text{such that}\  f(XY)=f(X)f(Y), \forall Y\in\A\ \text{or}\  f(YX)=f(X)f(Y), \forall Y\in\A\Big\}.
\end{equation*}
\endgroup
\end{corollary}
Given an element $X\in\A$, $A$-algebraic numerical range $V_A(X)$ of $X$ was defined in \cite{B.M.Positivity} by
\begin{align*}
V_A( X):= \big\{f(AX): f\in\mathcal{S}_{A}(\A)\big\}.
\end{align*}
Since $A$ is well-supported, then for any $X\in\A^A$ the $A$-numerical range is a nonempty, convex and compact subset of $\mathbb{C}$.
More details and properties of this set were presented in \cite{amz, B.M.Positivity, M.Z}.
As an immediate consequence of Theorems \ref{th2-2}, we have the following result.
\begin{corollary}
Let $\A$ be a \cs-algebra and $X\in\A^A$.
Then $\sigma_A(X)\subseteq V_A(X)$.
\end{corollary}
The next main result presents an $A$-spectral permanence property for \cs-algebras.
\begin{theorem}\label{perma}
Let $X\in\A^A$. Then
\begin{itemize}
\item[(i)] $\sigma_A^r(X)=\Big\{f(AX): f\in\mathcal{S}_{A}(\A), \,f(AXX^*A)=f(AX)f(AX^*A)=|f(AX)|^2f(A^2)\Big\}$.
\item[(ii)] $\sigma_A^l(X)=\Big\{f(AX): f\in\mathcal{S}_{A}(\A),\, f(X^*AX)=|f(AX)|^2\Big\}$.
\end{itemize}
In particular if $\B$ is a \cs-algebra such that $A\in\B\subseteq\A$, then
$\sigma_A(X,\A)=\sigma_A(X,\B)$ for any $X\in\B^A$.
\end{theorem}
\begin{proof}
(i) Let $\zeta\in\sigma_A^r(X)$. So, by Theorem \ref{th2-2},
there exists $f\in\mathcal{S}_{A}(\A)$ such that $f(AXY)=f(AX)g(AY)$ for any $Y\in\A^A$ and $\zeta=f(AX)$.
Since $X^*A\in\A_A\subseteq\A^A$, we have
\begin{align*}
f(AXX^*A)&=f(AX)f(AX^*A)=f(AX)\overline{f(AXA)}
\\& =f(AX)\overline{f(AX)f(A^2)} = \left|f(AX)\right|^2f(A^2).
\end{align*}
Hence $\zeta\in\Big\{f(AX): f\in\mathcal{S}_{A}(\A), \,f(AXX^*A)=f(AX)f(AX^*A)=|f(AX)|^2f(A^2)\Big\}$.

For the converse assume that there exists $f\in\mathcal{S}_{A}(\A)$ such that $f(AXX^*A)=f(AX)f(AX^*A)=|f(AX)|^2f(A^2)$.
Set $\eta=f(AX)$.
For any $Y\in\A^A$, by the Cauchy--Schwarz inequality, we have
\begin{align}\label{C.S.22}
\big|f\left(A(X-\eta\uno)Y\right)\big|^2\leq f\big(A(X-\eta\uno)(X^*-\overline{\eta}\uno)A\big)f\left(Y^*Y\right).
\end{align}
Since
\begingroup\makeatletter\def\f@size{10}\check@mathfonts
\begin{align*}
f\big(A(X-\eta\uno)(X^*-\overline{\eta}\uno)A\big)=f(AXX^*A)-\eta f(AX^*A)-\overline{\eta f(AX^*A)}+|\eta|^2f(A^2) = 0,
\end{align*}
\endgroup
we deduce from \eqref{C.S.22} that $f\left(A(X-\eta\uno)Y\right)=0$ for any $Y\in\A^A$;
that is $f(AXY)=f(AX)f(AY)$ for any $Y\in\A^A$.
Hence, by Theorem \ref{th2-2}(i), we get $\eta\in\sigma_A^r(X)$.

(ii) Suppose that there exists $f\in\mathcal{S}_{A}(\A)$ such that $f(X^*AX)=|f(AX)|^2$.
Let $\lambda=f(AX)$. For any $Y\in\A^A$, by Theorem \ref{adjoint}(ii),
there exists $Z\in\A^A$ such that $A^{1/2}Y=Z^*A^{1/2}$.
Then, by the Cauchy--Schwarz inequality, we have
\begin{align*}
&\big|f(AYX)-f(AX)f(AY)\big|^2
\\&\qquad \qquad=\big|f\left(AY(X-\lambda\uno)\right)\big|^2
\\&\qquad \qquad=\left|f\left(A^{1/2}Z^*A^{1/2}(X-\lambda\uno)\right)\right|^2
\\& \qquad \qquad\leq f\left(A^{1/2}Z^*ZA^{1/2}\right)f\left((X^*-\overline{\lambda}\uno)A(X-\lambda\uno)\right)
\\& \qquad \qquad= f\left(A^{1/2}Z^*ZA^{1/2}\right)\Big(f(X^*AX) -\overline{\lambda}f(AX) - \lambda f(X^*A)+ |\lambda|^2f(A)\Big)
\\& \qquad \qquad= f\left(A^{1/2}Z^*ZA^{1/2}\right)\Big(|\lambda|^2 -|\lambda|^2 - |\lambda|^2+ |\lambda|^2\Big) =0.
\end{align*}
Thus $f(AYX)=f(AX)f(AY)$ for any $Y\in\A^A$ and then by Theorem \ref{th2-2}(ii) we obtain $\lambda\in\sigma_A^l(X)$.

Now, let $\lambda\in\sigma_A^l(X)$. Then $A^{1/2}\A^A(X-\lambda\uno)\subsetneq A^{1/2}\A^A$.
A similar reasoning as in the proof of Lemma \ref{lem-1} with $A^{1/2}$ instead of $A$,
we can show that there exists a positive linear functional $g$ on $\A$ such that $g(A^{1/2})=1$ and $g(A^{1/2}Y(X-\lambda\uno))=0$ for all $Y\in\A^A$.
In particular $\lambda=g(A^{1/2}X)$ and $g(A^{1/2}YX)=g(A^{1/2}Y)g(A^{1/2}X)$ for all $Y\in\A^A$.
As $g(A^{1/2})=1$ then $g(A)\neq 0$. Set $f=\frac{g}{g(A)}$ and note that $f$ belongs to $\mathcal{S}_{A}(\A)$.
Also
\begin{align*}
g(A)f(AX)=g(AX)=g(A^{1/2}A^{1/2}X)=g(A)g(A^{1/2}X)=g(A)\lambda.
\end{align*}
Hence $\lambda=g(A^{1/2}X)=f(AX)$.
Note that by Theorem \ref{adjoint}(ii) we have $A^{1/2}X=W^*A^{1/2}$ for some $W\in\A^A$ and so
\begin{align*}
f(X^*AX)&= f(X^*A^{1/2}A^{1/2}X)
\\&= f(A^{1/2}WA^{1/2}X)=\frac{g(A^{1/2} WA^{1/2}X)}{g(A)}
\\&= \frac{g(A^{1/2}W)g(AX)}{g(A)}=g(X^*A^{1/2})f(AX)
\\&=\overline{f(AX)}f(AX)=|f(AX)|^2.
\end{align*}
Therefore, $\lambda\in\Big\{f(AX): f\in\mathcal{S}_{A}(\A),\, f(X^*AX)=|f(AX)|^2\Big\}$
and the proof (ii) is thus complete.

Finally, let $\B$ be a \cs-algebra such that $A\in\B\subseteq\A$ and let $X\in\B^A$.
Then $AX, X^*AX$ are in $\B^A$. By Hahn--Banach theorem we infer that
$\sigma_A(X,\A)=\sigma_A(X,\B)$.
\end{proof}
As an immediate consequence of Theorem \ref{perma}, we get the following result.
\begin{corollary}\label{lem2-0}
Let $X\in\A$. Then $\lambda\in\sigma^r(X)$ (resp. $\lambda\in\sigma^l(X)$) if and only if there
exists $f\in\mathcal{S}(\A)$ such that $\lambda=f(X)$ and $f(XX^*)=|f(X)|^2$ (resp. $f(X^*X)=|f(X)|^2$).
\end{corollary}
\begin{proof}
In Theorem \ref{perma} put $A=\uno$.
\end{proof}
In the following result, we give another characterization of $\sigma_A(X)$.
\begin{theorem}\label{mr1}
Let $X\in\A^A$ and $P$ be the orthogonal projection onto the range of $A$.
Then $\sigma_A(X)\backslash\{0\}=\sigma(PX)\backslash\{0\}$.
In particular $\sigma_A(P)=\{1\}$.
\end{theorem}
\begin{proof}
First observe that by Theorem \ref{perma}, we have $\sigma_A(P)=\{1\}$, since $AP=A$.
Let $\lambda\in\sigma(PX)\backslash\{0\}$. By Corollary \ref{cor-2}, there exists a state $f\in\mathcal{S}(\A)$
such that $\lambda=f(PX)$ and $f(PXY)=f(PX)f(Y)$ for all $Y\in\A$
or $f(YPX)=f(PX)f(Y)$ for all $Y\in\A$.
Using the fact that $PXP=PX=P^2X$, $f(PX)=f(PX)f(P)$.
Since $\lambda\neq 0$, we get $f(P)=\|P\|=1$.
In particular $f(PY)=f(YP)=f(P)f(Y)$ for any $Y\in\A$.
Whence $\lambda=f(X)$ and $f(XY)=f(X)f(Y)$ for any $Y\in\A$.
If $\lambda\notin\sigma_A(X)\backslash\{0\}$, then $A(\lambda\uno-X)Y=A$ for some $Y\in\A^A$.
Accordingly $P(\lambda\uno-X)Y=P$. Hence $1=f(P)=f( P(\lambda\uno-X)Y)=f(X)-\lambda=0$ which is a contradiction.
Now, let $\lambda\in\sigma_A(X)\backslash\{0\}$. Then by Theorem \ref{th2-2}, $\lambda=g(AX)$
with $g\in\mathcal{S}_{A}(\A)$ and $g(AXY)=g(AX)g(AY)$, for any  $Y\in\A^A$ or $g(AYX)=g(AX)g(AY)$, for any  $Y\in\A^A$.
To fix ideas, assume that $g(AXY)=g(AX)g(AY)$, for any $Y\in\A^A$ and set $\vf(Z)=\frac{g(AZA)}{g(A^2)}$, $Z\in\A$.
Since $YA$ is in $\A^A$ for any $Y\in\A$, straightforward computation entails that $\lambda=\vf(PX)$
and $\vf(PXY)=\vf(PX)\vf(Y)$ for any $Y\in\A$. Hence, by Corollary \ref{cor-2}, $\lambda\in\sigma(PX)\backslash\{0\}$.
The case when $g(AYX)=g(AX)g(AY)$, for any $Y\in\A^A$ can be handled in a similar way.
\end{proof}
\begin{remark}
Note that, in general, it may happen that $\sigma_A(X)\neq\sigma(PX)$ for some $X\in\A^A$.
For instance, if the range of $A$ is not dense in $\h$, then $P\notin\{0, I\}$. In particular $\sigma(P)=\{0, 1\}$.
Using Theorem \ref{perma} and the fact that $AP=A$, it is easy to see that $\sigma_A(P)=\{1\}$.
\end{remark}
Following \cite[p.~9]{Mu}, if $X$ is an element of $\A^A$, its $A$-spectral radius is defined to be
\begin{align*}
r_A(X):=\sup\big\{|\lambda|: \lambda\in\sigma_A(X)\big\}.
\end{align*}
\begin{theorem}\label{per}
Let $X\in\A^A$. Then $\sigma_A(X)$ is a nonempty compact subset of $\mathbb{C}$ and
\begin{align*}
r_A(X)=\displaystyle\lim_{n\to\infty}{\|X^n\|}^{1/n}_{A}.
\end{align*}
\end{theorem}
\begin{proof}
It is clear that from Theorem \ref{perma} that $\sigma_A(X)$ is a nonempty compact subset of $\mathbb{C}$.
Let $P$ be the orthogonal projection onto the range of $A$. Since $P=(A^{1/2})^\dag A^{1/2}$, by Theorem \ref{mr1}, we have
\begin{align*}
r_A(X)&=r(PX)=r\left((A^{1/2})^\dag A^{1/2}X\right)
=r\left(A^{1/2}X(A^{1/2})^\dag\right)
\\&=\displaystyle\lim_{n\to\infty}\left\|\left(A^{1/2}X(A^{1/2})^\dag\right)^n\right\|^{1/n}
=\displaystyle\lim_{n\to\infty}\left\|A^{1/2}X^n(A^{1/2})^\dag\right\|^{1/n}.
\end{align*}
Using \cite[Proposition 2.2]{arias2008metric} and \cite[Theorem~3.5]{B.M.Positivity},
we infer that $r_A(X)=\displaystyle\lim_{n\to\infty}{\|X^n\|}^{1/n}_{A}$.
\end{proof}
The next result is a variant of Gelfand--Mazur theorem for the $A$-spectrum.
\begin{corollary}
If every nonzero element of $\A^A$ is $A$-invertible, then $\A^A = \mathbb{C}\uno$.
In particular $A=\alpha \uno$ for some $\alpha>0$ and $\A^A=\A$.
\end{corollary}
\begin{proof}
Let $X\in\A^A$. By Theorem \ref{per} there is a complex number $\lambda$ such that $\lambda\in\sigma_A(X)$.
Thus $\lambda\uno-X$ is not $A$-invertible and so $\lambda\uno-X=0$. Hence $X=\lambda\uno$.
\end{proof}
For $X\in\A^A$, let $\partial\sigma_A(X)$ denote the boundary of $\sigma_A(X)$.
Here is another consequence of Theorem \ref{per}.
\begin{corollary}\label{tb}
Let $X\in\A^A$. If $\lambda\in \partial\sigma_A(X)$, then there exists a sequence $(X_n)$ in $\A^A$ such that ${\|X_n\|}_{A}=1$ and
\begin{align*}
\displaystyle\lim_{n\to\infty}{\big\|X_n\left(\lambda \uno-X\right)\big\|}_{A}=\displaystyle\lim_{n\to\infty}{\big\|\left(\lambda\uno-X\right)X_n\big\|}_{A}=0.
\end{align*}
\end{corollary}
\begin{proof}
Let $\lambda\in \partial\sigma_A(X)$.
Since by Theorem \ref{per} $\sigma_A(X)$ is a nonempty compact subset of $\mathbb{C}$,
there exists a sequence $(\lambda_n)$ of complex numbers
converging to $\lambda$ with $\lambda_n\notin\sigma_A(X)$. For any $n\ge 1$ let $Y_n\in\A^A$ such that $A(\lambda_n\uno-X)Y_n=AY_n(\lambda_n\uno-X)=A$.
Set $X_n=\frac{Y_n}{{\|Y\|}_{A}}$. Then ${\|X_n\|}_{A}=1$ and easy computations entails that
$\displaystyle\lim_{n\to\infty}{\big\|X_n\left(\lambda \uno-X\right)\big\|}_{A}=
\displaystyle\lim_{n\to\infty}{\big\|\left(\lambda\uno-X\right)X_n\big\|}_{A}=0$.
\end{proof}
Note that if $X\in\A^A$, by Theorem \ref{per} the sequence $(X_n)_{n\ge 1}$ in Corollary \ref{tb} is also in $\A_A$.
Further, observe that if $X\in\A_A$, then
\begin{equation}\label{inega}
\|A^{1/2}X\|\le \|A^{1/2}\|\|X\|_A.
\end{equation}
Indeed, using Remark \ref{rem1} and the definition of ${\|\!\cdot\!\|}_{A}$ we see that
$f(X^*AX)\leq {\|X||}_{A} f(A)$ for any $f\in\mathcal{S}(\A)$.
Taking the supremum over all $f\in\mathcal{S}(\A)$ we get $\|A^{1/2}X\|^2=\|X^*AX\|\leq {\|X\|}^2_{A}\|A\|$.
In addition, by \cite[Lemma~2.12]{amz} there exists a positive constant $c$ such that
\begin{equation}\label{est}
{\|X\|}_{A}\leq c\|A^2X\|,
\end{equation}
for all $A$-self-adjoint operators $X\in\A_A$.
We use these facts in the proof of the next result.
\begin{theorem}\label{bound}
Let $X\in\A^A$. Then $\sigma_A^r(X)\subseteq \sigma(X)$ and
$\partial\sigma_A(X)\subseteq\sigma_A^r(X)\cap \sigma^l_A(X)$.
In particular $\sigma_A^l(X)$ and $\sigma_A^r(X)$ are  nonempty  sets.
\end{theorem}
\begin{proof}
Let $\lambda\in\sigma_A^r(X)$. By Theorem \ref{perma}(i), there exists $f\in\mathcal{S}_{A}(\A)$ so that $\lambda=f(AX)$ and $f(AXX^*A)=f(AX)f(AX^*A)=|f(AX)|^2f(A^2)$. Put $\varphi=\frac{f}{f(A^2)}$.
Since $\A A\subset \A^A$, straightforward computation entails that $\lambda=\varphi(X)$ and $\varphi(XY)=\varphi(X)\varphi(Y)$ for any $Y\in\A$. Hence by Corollary \ref{lem2-0} we infer that $\lambda\in\sigma(X)$.

Now, let $\lambda\in\partial\sigma_A(X)$. By Corollary \ref{tb} there
exists a sequence $(X_n)$ in $\A_A$ such that ${\|X_n\|}_{A}=1$ and
\begin{align*}
\displaystyle\lim_{n\to\infty}{\big\|X_n\left(\lambda \uno-X\right)\big\|}_{A}=\displaystyle\lim_{n\to\infty}{\big\|\left(\lambda\uno-X\right)X_n\big\|}_{A}=0.
\end{align*}
Note that $X_n\in\A^A$, and then there exists $Y_n\in\A_A$ such that $A X_n=Y_n^*A$.
If $\lambda\notin\sigma_A^r(X)$, then there exists  $Y\in\A^A$ such that $A(\lambda\uno-X)Y=A$.
Hence
\begingroup\makeatletter\def\f@size{10}\check@mathfonts
\begin{align*}
\|AX_n\|=\|Y_n^*A\|=\|Y_n^*A(\lambda\uno-X)Y\|=\|AX_n(\lambda\uno-X)Y\|\leq \|AX_n(\lambda\uno-X)\|\|Y\|.
\end{align*}
\endgroup
This together with inequalities \eqref{inega} and \eqref{est} imply that $\displaystyle\lim_{n\to\infty}\|X_n\|_A=0$. This is a contradiction since $\|X_n\|_A=1$ for any $n\ge 1$. Whence $\lambda\in\sigma_A^r(X)$.
Similarly we can prove that $\lambda\in\sigma_A^l(X)$. Thus $\partial\sigma_A(X)\subseteq\sigma_A^r(X)\cap \sigma^l_A(X)$.
\end{proof}
The next result gives another expression of the $A$-spectral radius for elements of $\A_A$.
\begin{corollary}\label{bcd}
Let $X\in\A_A$. Then
\begin{align*}
r_A(X)=\sup\big\{|\lambda|: \lambda\in\sigma_A^l(X)\big\}=\sup\big\{|\lambda|: \lambda\in\sigma_A^r(X)\big\}.
\end{align*}
\end{corollary}
\begin{proof}
Set $\alpha=\sup\big\{|\lambda|: \lambda\in\sigma_A^l(X)\big\}$ and
$\beta=\sup\big\{|\lambda|: \lambda\in\sigma_A^r(X)\big\}$.
It is clear that $\alpha\leq r_A(X)$ and $\beta\leq r_A(X)$.
Now, since $r_A(X)$ is attained on $\partial\sigma_A(X)$,
by  Theorem  \ref{bound} we get $\alpha\geq r_A(X)$ and $\beta\geq r_A(X)$.
 Hence $r_A(X) = \alpha = \beta$.
\end{proof}
We close the paper with the following result.
\begin{corollary}\label{inclusion}
Let $X\in\A_A$. Then $r_A(X)\leq r(X)$.
\end{corollary}
\begin{proof}
By Corollary \ref{bcd}, we have
\begin{align*}
r_A(X)=\sup\big\{|\lambda|: \lambda\in\sigma_A^r(X)\big\}\leq \sup\big\{|\lambda|: \lambda\in\sigma^r(X)\big\}=r(X).
\end{align*}
\end{proof}
\noindent
 \textbf{Declarations.} : The work presented here was carried out in collaboration between all authors.
All authors contributed equally and significantly in writing this article. All authors have contributed
to the manuscript. All authors have read and agreed to the published version of the manuscript.
\\
\textbf{Competing interests.} The authors declare that they have no conflict of interest.
\bibliographystyle{amsplain}

\end{document}